\documentclass[12pt,notitlepage,twoside]{article}
\usepackage{latexsym}
\usepackage{amsmath}
\usepackage{amsfonts}
\usepackage{amssymb}
\renewcommand{\a }{\alpha }
\renewcommand{\b }{\beta }
\renewcommand{\d}{\delta }

\newcommand{\D }{\Delta }

\newcommand{\e }{\varepsilon }

\renewcommand{\l }{\lambda }

\newcommand{\n }{D }

\newcommand{\Sig }{\Sigma}

\renewcommand{\th }{\theta }

\newcommand{\ov}{\overline}

\newcommand{\be}{\begin{equation}}
\newcommand{\ee}{\end{equation}}
\newenvironment{pf}{\noindent{\bf Proof.}\enspace}{%\rule{2mm}{2mm}
\hfill$\Box$\medskip}

\newenvironment{pfn}[1]{\noindent{\bf Proof of {#1}\enspace}}{%\rule{2mm}{2mm}
\hfill$\Box$\medskip}
\newcommand{\R}{\mathbb{R}}
\newcommand{\Z}{\mathbb{Z}}

\newcommand{\N}{\mathbb{N}}

\newtheorem{thm}{Theorem}[section]
\newtheorem{pro}[thm]{Proposition}
\newtheorem{lem}[thm]{Lemma}

\newtheorem{cor}[thm]{Corollary}
\newtheorem{df}[thm]{Definition}
\numberwithin{equation}{section}
\let \n = \noindent

\author{  Hichem Chtioui$^a$, Mohameden Ould Ahmedou$^b$\footnote{ Corresponding author, { \tt ahmedou@analysis.mathematik.uni-tuebingen.de}} and Ridha Yacoub$^c$\\
{\footnotesize a: D{\'e}partement de Math{\'e}matiques,
Facult{\'e} des Sciences de Sfax, Route de Soukra, Sfax, Tunisie.}\\
{\footnotesize b: Mathematisches Institut der Universit\"{a}t
T\"{u}bingen,
 Auf der Morgenstelle 10, 72076
T\"{u}bingen, Germany.}\\ {\footnotesize c: Institut
pr\'eparatoire aux Etudes d'Ing\'enieur  de Monastir, Avenue Ibn
Al Jazzar, 5019 Monastir, Tunisie.}}

 \title { \Large \textbf{Existence and Multiplicity  results for
 the  prescribed   Webster Scalar Curvature Problem
 on three $C  R $ manifolds}}

\begin{document}

\date{ }

\maketitle

{\footnotesize \noindent {\bf Abstract.}  This paper is devoted to  the existence
 of contact forms
of prescribed Webster scalar curvature  on a $3-$dimensional  CR
compact manifold locally conformally CR equivalent to the unit
sphere $\mathbb{S}^{3}$ of $\mathbb{C}^{2}$. Due to Kazdan-Warner
type obstructions, conditions on the function $H$ to be realized
as a Webster scalar curvature have to be given. We prove new
existence results based on a new type of Euler-Hopf type formula.
Our argument gives an upper bound on the Morse index of the
obtained solution. We also give a lower bound on the number of
conformal contact forms
 having the same Webster scalar curvature. }\\
 {\footnotesize \noindent{\bf Mathematics Subject Classification
(2000) :}\quad 53C15, 53C21, 35J65, 18G35.}

%\tableofcontents

\section{Introduction }
\mbox{} Let $(M, \th)$ be  a strictly pseudoconvex  CR compact manifold of dimension
$2n + 1$ locally  CR equivalent to the unit sphere
$\mathbb{S}^{2n + 1}$ of $\,\mathbb{C}^{n + 1}$ with a contact form
$\th$, and let $K : M \to \R$ be a $C^3$ positive function. The
prescribed Webster scalar curvature on $M$ is to find suitable
conditions on $K$ such that $K$ is the Webster scalar curvature
for some contact form $\tilde\th$ on
 $M$,  CR equivalent to $\th$. If we set
$\tilde\th=u^{\frac{2}{n} }\th$, where $u$ is a smooth positive
function on $M$, then the above problem is equivalent to solving
the following equation $$ (P)\quad \begin{cases}
 L_{\th}u = \frac{n}{2(n+1)} K u^{1+\frac{2}{n}} \quad\mbox{ in }M \\
u>0 \qquad\mbox{ in } M,
\end{cases}
  $$
 where $$L_{\th}u= \D _{\th}u+ \frac{n}{2(n+1)}\mathcal{R}_{\th} u $$
$\D_{\th}$ is the sublaplacian operator on $(M,\th )$  and
$\mathcal{R}_{\th}$ is the Webster scalar curvature of  $(M,\th
)$. \\
Problem $(P)$ is the analogue of the prescribed  scalar curvature
problem on Riemannian manifolds. While the scalar curvature
problem in the Riemannian framework was  extensively studied (see
for example the monograph \cite{A} and the references therein),
only few results  were established for problem $(P)$ (see
\cite{FU}, \cite{G}, and \cite{MU}). On the contrary, the Yamabe
problem on CR manifolds, that is when $K$ is assumed to be
constant, was widely studied by various authors (see \cite{JL1},
\cite{JL2}, \cite{JL3}, \cite{G1} and \cite{GY}). \\
The  problem $(P)$ has a variational structure, however
the associated  Euler functional   does not satisfy the
Palais-Smale condition, that is, there exist noncompact sequences
along which the functional is bounded and its gradient goes to
zero.
 Moreover, there are topological obstructions  of Kazdan-Warner
 condition type
 to solve $(P)$, see \cite{GL}.
 Hence one does not expect to solve problem $(P)$ for all functions $K$, and
 so it is natural to ask:
 under which conditions on $K$  does $(P)$   have a  solution? \\
In \cite{MU}, Malchiodi and Uguzzoni considered the case where
$M=\mathbb{S}^{2n+1}$ the unit sphere of $\mathbb{C}^{n+1}$ and
gave a perturbative result for problem $(P)$, that is when $K$ is
assumed to be a small perturbation of a constant (see also
\cite{FU}). Their approach uses a perturbation method due
 to Ambrosetti and Badiale \cite{AmB}.
In \cite{G}, N. Gamara noticed , in analogy with the
$4$-dimensional Riemannian case, that there is a balance
phenomenon between the self interactions and the mutual
interactions of the functions failing to satisfy Palais-Smale
condition in the $3$-dimensional CR case   ( see \cite{BC'} and
\cite{BCCH} for the Riemannian case). In \cite{G}  the case where
$M$ is locally conformally CR equivalent to the CR Sphere of
$\mathbb{C}^2$   was considered (thus when $n=1$), and a
Euler-Hopf type criterion for $K$ was provided to find solutions
for $(P)$ . The method used in \cite{G} is due to Bahri and Coron
\cite{BC'}. It consists of studying the critical points at
infinity of the associated variational problem, computing their
total Morse index, and comparing this total index to the
Euler-Poincar\'e characteristic of the space of variations. \\
In this paper we revisit the three dimensional case, namely the
following equation
$$ (P_{K})\quad \begin{cases}
 L_{\th}u = \frac{1}{4} K u^{3} \quad\mbox{ in }M \\
u>0 \qquad\mbox{ in } M.
\end{cases}
  $$
Our goal here is to give new existence results which generalize the one
obtained by N. Gamara \cite{G} and also to give, in generic cases, a lower
bound of the number of contact forms of prescribed Webster-Tanaka scalar curvature $K$. \\
To state our results, we set the following notations. Let $G(a,.)$
be the Green's function of $L_\th$ on $M$ and $A_a$ the value of
the regular part of $G$  at $a$. Let $\mathcal{K}$ the set of
critical points of $K.$ We say  that
  $K$  satisfies the condition $(C_0)$ if it has  only nondegenerate critical points
such that $$ \frac{-\D_\theta K(y)}{3K(y)}
-2A_{y} \neq 0 \, \quad \forall  y \in \mathcal{K} $$
 Now, we introduce the
following set
\begin{equation}\label{e:i}
\mathcal{K}_+ =\{ y \in \mathcal{K} \, ; \,\frac{-\D_\th
K(y)}{3K(y)}   -2A_{y}>0  \, \}.
\end{equation}
For $p\in\N^*$ and for any $p$-tuple $\tau_p=(y_{1},...,y_{p}) \in
(\mathcal{K}_+)^p$ such that $y_{i}\neq y_{j}$ if $i\neq j$, we
define a matrix $M(\tau_p)=(M_{ij})_{1\leq i,j\leq p}$, by
\begin{equation}\label{f:Matrix}
 M_{ii}=\frac{-\D_\th K(y_{i})}{3K(y_{i})^2}-2\frac{
A_{y_{i}}}{K(y_{i})}, \qquad M_{ij}
=-\frac{2G(y_{i},y_{j})}{\left(K(y_{i})K(y_{j})\right)^{1/2}}
\quad \mbox{for } i\neq j.
\end{equation}
We denote by $\rho(\tau_p)$ the least eigenvalue of $M(\tau_p)$ and we
say that a function $K$ satisfies the condition ${\bf (C_1)}$ if
 for every $\tau_p \in (\mathcal{K}_+)^p$, we have that $\rho(\tau_p) \, \ne 0$.\\
We set
\begin{equation}\label{f:finfty}
\mathcal{F_1} := \{ \tau_p = (y_1,\cdots,y_p) \in (\mathcal{K}_+)^p \, ; \, \, \rho(\tau_p) \, > 0 \}
\end{equation}
and define an index
$
\iota: \mathcal{F_1} \to \Z
$
defined by
$$\iota(\tau_p) := p - 1 \, + \, \sum_{i=1}^p (3 - m(K,y_i)),
$$
where $m(K,y_i)$ denotes the Morse index of $K$ at its critical point $y_i$.\\
Now we state our main result.

\begin{thm}\label{t:main}

Let $0 < K \in C^2(M)$ be a  positive function satisfying the conditions $(C_0)$ and $(C_1)$.\\
If there exists $k \in \N$ such that
\begin{enumerate}
\item
$$
\sum_{\tau_p \in \mathcal{F_1} ; \iota(\tau_p) \leq k - 1 } (-1)^{\iota(\tau_p)} \, \neq \, 1,
$$
\item
$$
\forall \tau_p \in \mathcal{F_1}, \,  \iota(\tau_p) \neq k.
$$
\end{enumerate}
Then, there exists a solution $w$ to   the problem $(P_K)$ such
that:
$$
morse(w) \, \leq k,
$$
where $morse(w)$ is the Morse index of $w$, defined as the
dimension of the space of negativity of the linearized operator:
$$
\mathcal{L}_w(\varphi) := L_{\th} (\varphi) \, - 3 w^2 \varphi.
$$
Moreover, for generic $K$ it holds
$$
\# \mathcal{N}_k \, \geq | 1 \, - \, \sum_{\tau_p \in \mathcal{F_1} ; \iota(\tau_p) \leq k - 1 } (-1)^{\iota(\tau_p)} |,
$$
where  $\mathcal{N}_k$ denotes the set of solutions of $(P_K)$
 having their Morse indices less than or equal to $k$.
\end{thm}
Please observe that, taking  in the above $k$ to be $l_{\#} + 1,$
where $l_{\#} $ is the maximal index over all elements of
$\mathcal{F_1}$, the second assumption is trivially satisfied.
Therefore in this case, we have the following corollary, which
recovers previous existence results for three dimensional CR
manifolds locally CR equivalent to the sphere $\mathbb{S}^3$ of
$\mathbb{C}^2$ due to N. Gamara \cite{G}.

\begin{cor}\label{c:cor1}
Let $0 < K \in C^2(M)$ be a  positive function satisfying
the conditions $(C_0)$ and $(C_1)$.\\
If

$$
\sum_{\tau_p \in \mathcal{F_1} } (-1)^{\iota(\tau_p)} \, \neq \, 1,
$$
Then the problem $(P_K)$ has at least one solution.\\
Moreover, for generic $K$ it holds
$$
\# \mathcal{S} \, \geq | 1 \, - \, \sum_{\tau_p \in \mathcal{F_1}  } (-1)^{\iota(\tau_p)} |,
$$
where  $\mathcal{S}$ denotes the set of solutions of $(P_K)$.
\end{cor}

We point out that the main new contribution of Theorem
\ref{t:main} is that we
 address here  the case where the total sum in the above corollary  equals 1, but
 a partial one is not equal 1. The main issue being the possibility  to use such an
 information to prove existence of solution to the problem $(P_K)$.
Please notice that an  interpretation of the fact that the above
sum is different from one, is that the topological contribution of
the {\it critical points at infinity } to the topology of the
level sets of the associated Euler-Lagrange functional is not
trivial. In view of such an interpretation, we rise the following
question: what happens if the total contribution
 is trivial, but some  critical points at infinity induce a nontrivial difference of
 topology. Can we still use such a topological information to prove existence of solution ? \\
With respect to the above question,   theorem \ref{t:main} gives a sufficient
 condition to be able to derive from such a local information,  an existence as well as
 a multiplicity result together with information on the Morse index of the
 obtained  solution. In Section 5, we give a more general
 condition.
 Since this  condition involves the critical points at infinity of the variational problem,
 we have postponed its statement to the end of this paper, see please Theorem \eqref{t:a}. \\
As pointed out above, our result does not  only give existence results, but also,
under generic conditions, gives  a lower bound on the number of solutions of $(P_K)$.
Such a result is reminiscent to the celebrated Morse  Theorem, which states that, the
number of critical points of a Morse function defined on a compact manifold, is
lower bounded in terms of the topology of the underlying manifold. Our result can be seen
as some sort of {\it Morse Inequality at Infinity}. Indeed it gives a lower bound
on the number of metrics with prescribed curvature in terms of the
{\it topology at infinity}.\\
The remainder of this paper is organized as follows. In section 2
we set up the variational problem, its critical points at infinity
are characterized in Section 3. Section 4 is devoted to the proof
of the main result theorem \ref{t:main} while we give in Section 5
a more general statement than theorem \ref{t:main}.

\bigskip

\section{Variational setting and lack of compactness }
\mbox{} In this section we recall the functional setting and the
variational problem associated to $(P_K)$.\\
 Problem $(P_K)$ has a variational
structure, the functional being $$ J(u)= \frac {\int_{M}Lu\, u
\,\th \wedge d\th} {\left(\int_{M}K\,u^{4}\,\th \wedge d\th
\right)^{\frac{1}{2}}}, $$ defined on the unit sphere of
$\mathcal{S}^2_1(M)$ equipped with the norm
\begin{eqnarray}\label{e:01}
||u||^2=\int_{M}u \, Lu  \, \th \wedge d\th,
\end{eqnarray}
where $\mathcal{S}^2_1(M)$ is the Folland-Stein space (see
\cite{FS} for  definition).\\
Problem $(P_K)$ is equivalent to finding the critical points of
$J$ subjected to the constraint $u\in \Sig^+$, where
\begin{eqnarray}\label{e:02}
\Sigma ^+=\{u \in \Sigma \, / \, u\geq 0\}, \quad \Sig=\{u\in
\mathcal{S}^2_1(M)/\, \, ||u||=1\}
\end{eqnarray}
The Palais-Smale condition fails to be satisfied for $J$ on
$\Sigma^+$. To characterize  the sequences failing the
Palais-Smale condition, we need to set  some notations and
constructions.\\
Since $M$ is compact and locally CR
equivalent to $\mathbb{S}^{3}$, any point $a$ in $M$ has a
neighborhood $U_a\supset B_r(a)$, $r$ is independent of $a$, where
CR normal coordinates are defined, and such that the contact form
of $M$ is conformal to the standard contact form $\th_0$ of the
Heisenberg group $\mathbb{H}^1$; that is there exists a positive
function $\tilde{u}_a$ on $B_r(a)$ such that
$\th_0=\tilde{u}_a^{2}\th$, ($\tilde{u}_a$ smoothly
dependent on $a$). Let $u_a(x)=w_a(x)\tilde{u}_a(x)$, where
$w_a(x)=\chi(|x|)$, $\chi$ is a cut-off function $ \chi: \R
\rightarrow [0,1]$ defined by $$ \chi (t)= 1\quad \mbox{if}\quad
0\leq t\leq r/2;\quad \chi (t)= 0\quad \mbox{if}\quad t\geq r $$
and $\;|x|=|\exp_a^{-1}(x)|_{\mathbb{H}^1}$, where, letting
$(z,t)=\mbox{exp}_a^{-1}(x)$, $\mbox{exp}_a$ being the parabolic
exponential map based at $a$, then
$|(z,t)|_{\mathbb{H}^n}=(|z|^4+t^2)^{\frac{1}{4}}$ is the norm of
the Heisenberg  group  $\mathbb{H}^1$ ( see please \cite{JL1},
\cite{JL2}).\\

\noindent
Let $\l$ be a large positive parameter. We introduce on $B_r(a)$
the function
\be\label{e:o3}
 \d_{(a,\l)}(x)=c_1 \l|1+\l^2(|z|^2-it)|^{-1},
 \ee  and the constant $c_1$ is chosen such that the following
equation is satisfied $$ L_{\th_0}\d_{(a,\l)}=
\d_{(a,\l)}^{3}\quad\mbox{on}\quad B_r(a). $$
 Let
\begin{equation}\label{delta'}
\hat{\d}_{(a,\l) }(x)=\begin{cases}
 u_a\d_{(a,\l)}(x)\quad \mbox{in}\quad
B_r(a)\\ 0\quad \mbox{in}\quad B_r(a)^c.
\end{cases}
\end{equation}
 We  define a family of "almost solutions" $\tilde\d_{(a,\l )}$
to be the unique solution of
 \begin{equation*}\label{a} L_{\th}\tilde{\d}_{(a,\l) }(x)
=\bigl(\hat{\d}_{(a,\l)}(x)\,\bigr)^{3} \quad \mbox{ in } M.
\end{equation*}

\noindent
Setting
$$
H_{a,\l} := \l (\tilde{\d}_{a,\l} - \hat{\d}_{a,\l}),
$$
we have that:
\begin{pro}\label{p:h} $\cite{G}$
For $\l$ large, there exists a constant $C =C(\varrho)$ such that:
$$
|H_{a,\l}|_{L^{\infty}} \, \leq C ; \quad  \l |\frac{\partial H_{a,\l}}{\partial \l}|_{L^{\infty}} \leq C ; \quad  \l^{-1} |\frac{\partial H_{a,\l}}{\partial a}|_{L^{\infty}} \leq C .
$$
Moreover for $\varrho$ small and $\l$ large there holds:

\begin{eqnarray}
H_{a,\l}(a) \to A_a \, \, \mbox{ as }  \l \to \infty \\
H_{a,\l}(x) \to G(a,x) \, \, \mbox{ outside } B_{2 \varrho}(a) \, \mbox{ as }  \l \to \infty,
\end{eqnarray}
where $G(a,x)$ is the Green's function of the conformal
subLaplacian $L_{\th}$ and $A_a$ the value of its regular part
evaluated at $a.$
\end{pro}

\noindent
We define now the set of potential critical points at infinity
associated to the functional $J$.\\
For $\e>0$ and $p\in \N^*$, let us define
\begin{align*}
V(p,\e )=& \Big\{u\in \Sig /\exists a_1,...,a_p \in M, \exists \l
_1,...,\l _p >0, \exists \a _1,...,\a _p>0
 \mbox{ s.t. }\|u-\sum_{i=1}^p\alpha_i\tilde{\d}_{(a_i,\l_i)}\|<\varepsilon,
 \\& |\frac{\alpha_i^{2}K(a_i)}{\alpha_j^{2}K(a_j)}-1|<\e,
 \quad  \e_{ij}<\e \quad \l_i >\e ^{-1}\Big\},
\end{align*}
where $\e _{ij}^{-1}=\Big(\l _i/\l _j +\l _j/\l _i + \l _i\l
_jd(a_i,a_j)^2\Big)$, and  $
d(x,y)=\displaystyle{|exp_x^{-1}(y)|_{\mathbb{H}^1}}$ if $x$ and
$y$ are
 in a small ball of $M$ of radius $r$,  and $ d(x,y)$ is equal to
   ${r\over 2}$  otherwise.\\
\noindent For $w$ a  solution of $(P_K)$ we also define
$V(p,\e,w)$ as
\begin{eqnarray}
\{u\in \Sig/\exists \, \a_0>0 \mbox { s. t. } u-\a_0w\in
V(p,\e)\mbox{ and } |\a_0^2
J(u)^2 \, -1|<\e\}.
\end{eqnarray}
\noindent
The failure of Palais-Smale condition can be
described, following the ideas introduced in \cite{BrC} \cite{L}
 \cite{S}, as follows:
\begin{pro}\label{p:palais-smale}
 Let $(u_j)\in \Sig^+$ be a sequence such that $\n
J'(u_j)$ tends to zero and $J(u_j)$ is bounded. Then, there exist
an integer $p\in \N^*$, a sequence $\e_j>0$, $\e_j$ tends to zero,
and an extracted subsequence of $u_j$'s, again denoted $u_j$, such
that $u_j\in V(p,\e_j,w)$ where $w$ is zero or a solution of
$(P_K).$
\end{pro}

If a function $u$ belongs to $V(p,\e)$, we consider the following
minimization problem for $u\in V(p,\e)$ with $\e$ small
\begin{eqnarray}\label{e:51}
\min\{||u-\sum_{i=1}^p\a _i\tilde{\d}_{(a_i,\l_i)} ||,\, \a
_i>0,\, \l _i>0,\, a_i\in  M\}.
\end{eqnarray}

We then have the following proposition which defines a
parameterization of the set $V(p,\e )$. It follows from
corresponding  statements in \cite{B2}, \cite{BC}.

\begin{pro}\label{p:25}
For any $p\in \N^*$, there is $\e _p>0$ such that if $\e <\e _p$
and $u\in V(p,\e )$, the minimization problem \eqref{e:51} has a
unique solution (up to permutation). In particular, we can write
$u\in V(p,\e )$ as follows $$ u=\sum_{i=1}^p\bar{\a
}_i\tilde{\d}_{\bar{a}_i,\bar{\l }_i}+ v, $$ where $(\bar{\a
}_1,...,\bar{\a }_p,\bar{a}_1,...,\bar{a}_p,\bar{\l }_1,...,
\bar{\l }_p)$ is the solution of \eqref{e:51} and $v\in
\mathcal{S}^2_1(M)$ such that
 $$(V_0)\qquad < v,\psi>_\th=0 \ \mbox{ for
all }\  \psi\in
\bigg\{\tilde\d_i,\frac{\partial\tilde\d_i}{\partial\l_i},\frac{\partial\tilde\d_i}{\partial
a_i}, \text{ for } i=1, \ldots, p\bigg\}. $$
\end{pro}
Here, $<\  ,\  >_{\th}$ denotes the $L_\th$-scalar product defined
on $\mathcal{S}^2_1(M)$ by
\begin{eqnarray}\label{e:2.6}
< u,\, v>_{\th} = \int_M L_\th u\,v\,\th \wedge d\th.
\end{eqnarray}
Let $\nabla\!_\th$ be the CR gradient (or  subelliptic gradient)
which can be characterized  by
\begin{eqnarray}\label{e:2.7}
 \int_M \nabla\!_\th u \,\nabla\!_\th v\, \th\! \wedge\! d\th
=\int_M \Delta_\th u\,v\,\th \wedge d\th.
\end{eqnarray}

\noindent
In the following we will say that $v\in (V_0)$ if $v$ satisfies
$(V_0)$.
\begin{pro}\label{p:24} $\cite{G}$
There exists a $C^1$ map which, to each $(\a_1,..., \a_p, a_1,...,
a_p, \l_1,..., \l_p)$ such that $ \sum_{i=1}^p\a _i\tilde{\d}
_{(a_i,\l_i)} \in V(p,\e )$ with small $\e $, associates
$\ov{v}=\ov{v}_{(\a_i,a_i,\l_i )}$ satisfying
$$
J\left(\sum_{i=1}^p\a _i\tilde{\d} _{(a_i,\l_i)} +\ov{v}\right)= \min_{ v
\in (V_0)} J\left( \sum_{i=1}^p\a _i\tilde{\d} _{(a_i,\l_i)}
 +v\right).
$$
Moreover, there exists $c>0 $ such that the following holds
$$
||\ov{v}||\leq c \left(\sum_{i\leq p}(\frac{|\nabla\!_\th
K(a_i)|}{\l_i}+\frac{1}{\l_i^2})+\sum_{k\ne r}\e _{kr} (Log (\e
_{kr}^{-1}))^{1/2}\right).
$$
\end{pro}

\noindent Let $w$ be a  solution of $(P_K)$. The following
proposition defines a parameterization of the set $V(p,\e,w)$. Its
proof follows from the same arguments used to prove similar
statements in \cite{B2}.

\begin{pro}\label{p:w-parametrisierung}
There is $\e_0>0$ such that if $\e\leq \e_0$ and $u\in V(p,\e,w)$,
then the problem
 $$
 \min_{\a_i>0,\,,\,  \l_i>0,\, \,  a_i\in  {M} , \, \, h\in
T_w(W_u(w))}  \big|\big|u-\sum_{i=1}^p\a_i
\tilde{\d}_{(a_i,\l_i)}-\a_0(w+h) \big|\big|
 $$
has a unique solution $(\ov{\a},\ov{\l},\ov{a}, \ov{h} )$. Thus,
we write $u$ as follows:
$$u=\sum_{i=1}^p\ov{\a}_i\tilde{\d}_{(\ov{a}_i,\ov{\l}_i)}
+\ov{\a}_0(w+\ov{h})+v,$$ where $v$ belongs to
$H^1(M)\cap T_w(W_s(w))$ and it satisfies $(V_0)$,
$T_w(W_u(w))$ and $T_w(W_s(w))$ are the tangent spaces at $w$ to
the unstable and stable manifolds of $w$.
\end{pro}

\section{Critical points at infinity of the variational problem}

In the sequel, $\partial J$ designates the gradient of $J$ with
respect to the $L_{\th}$-scalar product $<\  ,\
>_{\th}$, that is $\forall\,u,\,v\in\mathcal{S}^2_1(M)$, we have
$<\partial J(u),\,v>_{\th}=J'(u)\,v.$ \\

Following A. Bahri we set the following definitions and
notations
\begin{df}
{\it A critical point at infinity} of $J$ on $\Sig^+$ is a limit
of a flow line $u(s)$ of the equation:
$$
\begin{cases}
\frac{\partial u}{\partial s} = - \n \partial J(u) \\
u(0) = u_0
\end{cases}
$$
such that $u(s)$ remains in $V(p,\e(s),w)$ for $s \geq s_0$.\\
Here $w$ is either zero or a solution of $(P_K)$ and $\e(s)$ is
some function tending to zero when $s \to \infty$. Using
Proposition \ref{p:w-parametrisierung}, $u(s)$ can be written as:
$$
u(s) \, = \,  \sum_{i=1}^p\a_i(s) \, \tilde{\d}_{(a_i(s),\l_i(s))} +
\a_0(s)(w + h(s)) \, + v(s).
$$
Denoting $a_i := \lim_{s \to \infty} a_i(s)$ and $\a_i=\lim_{s\to
\infty}\a_i(s)$, we denote by
$$(a_1,\cdots,a_p,w)_\infty \  \mbox{ or }\   \sum_{i=1}^p\a_i \,
\tilde{\d}_{(a_i,\infty)} + \a_0w$$
 such a critical point at infinity. If
$w \ne 0$ it is called of {\it $w$-type}.
\end{df}

\subsection{Ruling out the existence of critical points at infinity
in $V(p,\e,w)$ for $w \neq 0$}

The aim of this subsection is to prove  that, given a $C^2$
positive function $K$  satisfying the conditions of theorem
\ref{t:main} and  a solution $w$  of $(P_K)$, then for each $p \in
\N$, there are no critical point or critical point at infinity of
$J$ in the set $V(p,\e,w)$. The reason is that there exists a
pseudogradient of $J$ such that the Palais-Smale condition is
satisfied along
its decreasing flow lines.\\
In this section, for $u\in V(p,\e,w)$, using Proposition
\ref{p:w-parametrisierung}, we will write $u=\sum_{i=1}^p\a _i\,
\tilde{\d} _{(a_i,\l _i)}+\a_0(w+h)+v$.
\begin{pro}\label{p:32}
For $\e >0$ small enough and $u=\sum_{i=1}^p\a _i\tilde{\d} _{(a_i,\l
_i)}+\a_0(w+h)+v\in V(p,\e,w )$, we have the following expansion
\begin{align*}
J(u) = & \frac{S\sum_{i=1}^p \a_i^2+\a_0^2||w||^2}{ (S\sum_{i=1}^p
\a _i ^4K(a_i)+\a_0^4||w||^2 )^\frac{1}{2}}
\left[1-\frac{c_2\a_0}{\gamma_1}\sum_{i=1}^p
\a_i\frac{w(a_i)}{\l_i}\right.\\
& -\frac{1}{\gamma_1}\sum_{i\ne j} \a _i\a _j c_{ij}\e _{ij}
 +f_1(v)+Q_1(v,v) +f_2(h) +\a_0^2Q_2(h,h)\\
& \left.+o\biggl(\sum_{i\ne j}\e _{ij} +\sum_{i=1}^p\frac{1} {\l
_i }+||v||^2+ ||h||^2\biggr)\right]
\end{align*}
where
\begin{align*}
Q_1(v,v) = & \frac{1}{\gamma_1}||v||^2-\frac{3}{\b_1}
\int_{M}K\left(\sum_{i=1}^p(\a _i\tilde{\d}
_i)^2+(\a_0w)^2\right)
v^2, \\
Q_2(h,h)= & \frac{1}{\gamma_1}||h||^2-\frac{3}{\b_1}
\int_{M}K(\a_0w)^2 h^2,\\
f_1(v) = & -\frac{1}{\b_1 }\int_{M}K(\sum_{i=1}^p\a
_i\tilde{\d} _i)^3 \, v,  \\
f_2(h)= &
\frac{\a_0}{\gamma_1}\sum_i\a_i(\tilde{\d}_i,h)-\frac{\a_0}{\b_1}
\int_{M}K(\sum_i\a_i\tilde{\d}_i+\a_0 w)^3h,\\
c_2= & c_1^3
\int_{\mathbb{H}^1}\frac{1}{|1+|z|^2-it|^{3}}\th_0\wedge  d\th_0 ,
\qquad S = c_1^4
\int_{\mathbb{H}^1}\frac{1}{|1+|z|^2-it|^{4}}\th_0\wedge  d\th_0\\
 \b_1  = & {S}(\sum_{i=1}^p
\a _i^4 K(a_i))+\a_0^4||w||^2, \qquad
\gamma_1 = {S}(\sum_{i=1}^p \a _i^2)+\a_0^2||w||^2,
\end{align*}
and where $c_{ij}$ are $>0$ bounded constants.
\end{pro}
\begin{pf}
To prove the proposition, we
need to estimate
$$
N(u)= ||u||^2 \quad \mbox{ and }\quad D^2=\int_{M} K(x) u^4 \th\wedge  d\th,
$$
where $$ ||u|| := \int_M u L_{\th} u \, \th \wedge d \th$$

Now expanding $N(u)$, we get
$$
N(u):= \sum_{i=1}^p\a_i^2||\tilde{\d}_i||^2+2\a_i\a_0<\tilde{\d} _i,
w+h>_{\th} +\a_0^2 (||h||^2+||w||^2) +||v||^2 +\sum_{i\ne
j}\a_i\a_j<\tilde{\d} _i,\tilde{\d}_j>_{\th}.
$$
Now it follows from \cite{G} and elementary computations that
\begin{align}
& ||\tilde{\d}_i ||^2= S \, + c_2
\frac{H_{a_i,\l_i}(a_i)}{\l_i^2} \, + \, o(\frac{1}{\l_i^2}) ; \\
& <\tilde{\d}_i,\tilde{\d}_j>_{\theta} \, = \,  c_2\frac{
H_{a_j,\l_j}(a_i)}{\l_i \, \l_j} \, +  \,
 c_{ij} \e _{ij}(1 + o(1)) , \quad \mbox{ for } i\ne j, \label{e:o1}\\
 & <\tilde{\d}_i,w>_{\th}=\int_{B_r(a_i)}w\,\hat{\d}_i ^{3}\th \wedge d \th =
c_2\frac{w(a_i)}{\l_i}+ o(\frac{1}{\l_i}).
\end{align}
 Therefore
\begin{align}\label{N}
 N=  & \gamma_1+2\a_0\sum_{i=1}^pc_2
\a_i\frac{w(a_i)}{\l_i}+\a_i<\tilde{\d}_i,h>_{\th} +\sum_{i\ne
j}\a
_i \a _j c_{ij}\e _{ij}\\
 &  +\a_0^2||h||^2+ ||v||^2+
o\biggl(\sum_{i=1}^p\frac{1}{\l_i}+ \sum_{i\ne j}\e
_{ij}\biggr).\notag
\end{align}

\noindent
Now concerning the denominator, we compute it as follows
\begin{align}
D^2=&\int K (\sum_{i=1}^p\a _i\tilde{\d} _i
)^4\th \wedge d \th +\int K
(\a_0w)^4\th \wedge d \th \\
& + 4 \a_0 \, \int K (\sum_{i=1}^p\a _i\tilde{\d} _i)^3 w \th
\wedge d \th\,  + \,  4 \a_0^3\int K
(\sum_{i=1}^p\a _i\tilde{\d} _i)w^3\th \wedge d \th \notag\\
&  + 4 \int K
(\sum_{i=1}^p\a _i\tilde{\d} _i+\a_0w)^3(\a_0h+v) \th \wedge d \th \notag\\
 & + 12\int K(\sum_{i=1}^p\a_i \tilde{\d} _i+\a_0w)^2
(\a_0^2h^2+v^2+2\a_0hv)\th \wedge d \th  \notag\\
 & +O(\sum_{i=1}^p \int
 w^2\a_i^2\tilde{\d}_i^2+\a_0^2w^2\tilde{\d}_i^2)+O(||v||^{3}+\|h\|^3).\notag
\end{align}
 Observe that
\begin{align}
\int_{M}K (\sum_{i=1}^p\a _i\tilde{\d} _i)^4 \th \wedge d \th = &
\sum_{i=1}^p\a_i ^4K(a_i) S\\
& + 4   \sum_{i\ne j} \a_i^3\a_j K(a_i)c_{ij}\e_{ij} +O(
\frac{1}{\l _i^2})+o(\e _{ij}), \notag
\end{align}

\begin{align}
& \int_{M}K w^4\th \wedge d \th =||w||^2; \quad \int_{M} Kw^3\,
\tilde{\d}_i\, \th \wedge d \th  =
c_2\frac{w(a_i)}{\l_i }+ o(\frac{1}{\l_i }) ,\\
& \int_{M} K(\sum \a_i\tilde{\d}_i) ^3\, w\,\th \wedge d \th
=c_2\sum \a_i^3
K(a_i)\frac{w(a_i)}{\l_i } \, + \, o(\frac{1}{\l_i}),\\
& \int_{M}w^2\a_i^2\tilde{\d}_i^2+\a_0^2w^2\tilde{\d}_i^2\, \th
\wedge d \th =o(\frac{1}{\l_i }),
\end{align}

\begin{align}
\int_{M} K(\sum \a_i \tilde{\d}_i  + \a_0 w)^2  \,  v \, h\, \th
\wedge d \th  & = O\left(\int
\left(\sum \tilde{\d}_i^2 \,  +w^{-1}\sum  \tilde{\d}_i\right)|v||h|\right)  \notag\\
& =O\left(\|v\|^3+\|h\|^3+1/\l_i^{3}\right),
\end{align}
where we have used that $v\in T_w(W_s(w))$ and $h$ belongs to
$T_w(W_u(w))$ which is a finite dimensional space. Hence it
implies that $\|h\|_\infty \leq c\|h\|$.\\
Concerning the linear form in  $v$, since $v\in T_w(W_s(w))$, it
can be written as
\begin{align}
\int_{M}K (\sum_{i=1}^p & \a _i \tilde{\d} _i +\a_0
w)^{3}v\th \wedge d \th  \notag\\ & =\int K (\sum_{i=1}^p\a _i\tilde{\d}
_i)^{3}v \th \wedge d \th
+O\biggl(\sum_{i=1}^p \int(\a_i^2\a_0\tilde{\d}_i^2 \, w +\a_0^2\a_i
\tilde{\d}_i w ^2 )|v|\biggr)\notag\\
 & =f_1(v)+O\left(\frac{||v||}{\l_i}\right).
\end{align}
Finally, we have
\begin{align}
 \int K(\sum_{i=1}^p\a _i\tilde{\d} _i+\a_0w)^2 \, h^2\th \wedge d \th  = &
\a_0^2 \int
Kw^2 \, h^2 + o(||h||^2)\th \wedge d \th \\
 \int K(\sum_{i=1}^p\a _i\tilde{\d}_i+\a_0w)^2 \, v^2\th \wedge d \th = &
\sum_{i=1}^p\int K(\a _i\tilde{\d}_i)^2 \, v^2+\a_0^2
\int Kw^2 v^2\notag\\
 & +o(||v||^2).\label{v2}
\end{align}
Combining \eqref{N} to \eqref{v2}, and the fact that
 $\frac{\a_i^2K(a_i)}{\a_j^2K(a_j)}=1+o(1)$, the result
follows.
\end{pf}

Now, we state  the following lemma whose proof follows the
arguments used to prove similar statements in \cite{B1}, see the
Appendix of \cite{G} were the necessary modifications are given.
\begin{lem}\label{p:quadform}$\cite{G}$
 We have \\
(a) $Q_1(v,v)$ is a quadratic form positive definite in\\
\centerline{ $E_v=\{v\in \mathcal{S}_1^2(M)/v\in T_w(W_s(w)) \mbox{
and } v \mbox{ satisfies } (V_0)\}$.}
 (b) $Q_2(h,h)$ is a quadratic form
negative definite in $T_w(W_u(w))$.
\end{lem}

\begin{cor}\label{c:}  Let $u=\sum_{i=1}^p\a _i\d _{(a_i,\l
_i)}+\a_0(w+h)+v\in V(p,\e,w )$. There is an optimal
$(\ov{v},\ov{h})$ and a change of variables $v-\ov{v}\to V$ and
$h-\ov{h}\to H$ such that
$$J(u) = J\left(\sum_{i=1}^p\a _i\d _{(a_i,\l
_i)}+\a_0 w+ \ov{h}+\ov{v}\right)+||V||^2-||H||^2.
$$
Furthermore we have the following estimates
$$
||\ov{h}|| \leq \sum_i\frac{c}{\l_i }\, \, \, \mbox{ and }\, \, \,
||\ov{v}||\leq c\sum_i\frac{|\nabla_{\!\!\th}
K(a_i)|}{\l_i}+\frac{c}{\l_i^2}+c \sum_{k\neq
r}\e_{kr}(Log\e_{kr}^{-1})^{\frac{1}{2}},
$$
\begin{align*}
J(u)= & \frac{S_n\sum_{i=1}^p\a_i^2+\a_0^2||w||^2}{
(S_n\sum_{i=1}^p\a _i ^4 \, K(a_i)+\a_0^4 ||w||^2 )^{\frac{1}{2}}}
\left[1-\frac{c_2\a_0}{\gamma_1}\sum_{i=1}^p\a_i
\frac{w(a_i)}{\l_i }\right.\\
& -\frac{1}{\gamma_1}\sum_{i\ne j}\a _i\a _jc_{ij}\e _{ij}
\left.+o\left(\sum_{i\ne j}\e _{i j}+  \sum_{i=1}^p\frac{1} {\l _i
}\right)\right] +||V||^2-||H||^2.
\end{align*}
\end{cor}
\begin{pf}
Now observe that the above expansion of $J$ with respect to $h$ (respectively to $v$) is
almost equal  to $Q_2(h,h)+f_2(h)$
(respectively $Q_1(v,v)+f_1(v)$). As $Q_2$ is negative definite
(respectively $Q_1$ is positive definite), there is a unique
maximum $\ov{h}$ in the space of $h$'s (respectively a unique
minimum $\ov{v}$ in the space of $v$). Moreover , it is not difficult  to
prove that  $||\ov{h}||\leq c||f_2||$ and $\|\ov{v}\|\leq c\|f_1\|$.
Therefore the  estimate of $\ov{v}$ follows from Proposition \ref{p:24} while for 
the estimate of $\ov{h}$, we use the fact that for each $h\in
T_w(W_u(w))$ which is a finite dimensional space, we have
$\|h\|_\infty \leq c \|h\|$. Hence it follows  that
$\|f_2\|=O(\sum \l_i^{-1})$ and  our result follows.
\end{pf}

\noindent
Now we state the following corollary, which follows
 immediately from the above corollary and
the fact that $w > 0$ in ${M}$.

\begin{cor}\label{c::}
Let $K$ be a $C^2$ positive function and let $w$ be a
nondegenerate critical point of $J$ in $\Sig^+$. Then, for
each $p\in \N^*$, there is no critical
points or critical points at infinity in the set $V(p,\e,w)$, that
means we can construct a pseudogradient of $J$ so that the
Palais-Smale condition is satisfied along the decreasing flow
lines.
\end{cor}

\noindent Now, once the existence of mixed critical points at
infinity is ruled out, it follows from \cite{G}, that the critical
points at infinity are in one to one correspondence with the
elements of the set $\mathcal{F_1}$ defined in \eqref{f:finfty}.
That is,  a critical point at infinity corresponds to $\tau_p:=
(y_1,\cdots,y_p) \in (\mathcal{K}_+)^p$ such that the related
Matrix $M(\tau_p)$ defined in \eqref{f:Matrix} is positive
definite. Such a critical point at infinity will be denoted by
$\tau_p^{\infty} := (y_1,\cdots,y_p)_{\infty}.$ \\
Like  a usual critical point, it is associated to
 a {\it critical point at infinity} $x_{\infty}$ of the problem $(P_K)$,
 (which is a combination of classical critical points of $K$ with a $1-$dimensional asymptote),
 stable and unstable
manifolds, $W_s^{\infty}(x_{\infty})$ and
$W_u^{\infty}(x_{\infty})$.   These manifolds can be easily
described once a Morse type reduction is performed, see \cite{G}.
The stable manifold is, as usual, defined to be the set of points
attracted by the asymptote, under the action of the flow (see
below). The unstable one is a shadow object, which is the limit of
$W_u(x_{\l})$, $x_{\l}$ being the critical point of the reduced
problem and $W_u(x_{\l})$ its associated unstable manifolds.
Indeed the flow in this case splits the variable $\l$ from the
other variables near $x_{\infty}$. Notice that the flow of which
it is question above is the flow of a {\it pseudogradient at
infinity of Morse-Smale type }  for $-J$. Such a {\it
pseudogradient at infinity}, whose  existence  is ensured by the
Proposition 4.3 in \cite{G}, is known to have a very nice behavior
around the
critical points at infinity.\\
 We then may define the Morse index  $morse(x_{\infty})$ of the critical point at infinity
 $x_{\infty}$ to be equal to the dimension of
$W_u^{\infty}(x_{\infty})$. Observe that we have:
$\,morse(\tau_p^{\infty})=\iota(\tau_p)$. \\
 In the following definition, we extend the notion of domination of
critical points to critical points at infinity.
\begin{df}
$z_\infty$ is said to be dominated by another critical point at
infinity $z'_\infty$ if
$$W_u(z'_\infty)\cap W_s(z_\infty)\ne \emptyset.$$
\end{df}
We then write $\,  z_\infty<z'_\infty$. If we assume that the
intersection is transverse, then we obtain
$$morse(z'_\infty)\geq morse(z_\infty)+1.$$

\section{Proof of the main result}

This section is devoted to the proof of the main result of this paper,
theorem \ref{t:main}.

\begin{pfn}{\bf Theorem \ref{t:main} }

Setting
$$
l_{\#} := \sup \{ \iota(\tau_p); \, \tau_p \in \mathcal{F_1}  \}
$$
For $l \in \{0,\cdots,l_{\#} \}$ we define the following sets:
\begin{equation}\label{f:xl}
X^{\infty}_l := \bigcup_{\tau_p \in \mathcal{F_1} ; \,
\iota(\tau_p) \leq l} \ov{W^{\infty}_u(\tau_p^{\infty})},
\end{equation}
where $W_u^{\infty}(\tau_p^{\infty})$ is the unstable manifold
associated to the critical point at infinity $\tau_p^{\infty}$,
and
\begin{equation}\label{f:cl}
C(X^{\infty}_l) := \{ t \, u \, + \, (1 -t) \, (y_0)_{\infty},\,
 t \in [0,1],\,  u \in    X^{\infty}_l \},
\end{equation}
where $y_0$ is  a global maximum of $K$ on the manifold $M$.\\
By a theorem of Bahri-Rabinowitz \cite{BR}, it follows that:

$$
\ov{W^{\infty}_u(\tau_p^{\infty})} \, = \,
W^{\infty}_u(\tau_p^{\infty}) \,
 \cup \,  \bigcup_{x_{\infty} < \tau_p^{\infty}} W^{\infty}_u(x_{\infty}) \, \cup \,
  \bigcup_{w < \tau_p^{\infty}} W_u(w),
$$
where $x_{\infty}$ is a critical point at infinity dominated by
$\tau_p^{\infty}$ and $w$ is a solution of $(P_K)$ dominated by
$\tau_p^{\infty}.$ By transversality arguments  we assume that the
Morse index of $x_{\infty}$ and the Morse index of $w$ are not
bigger than $l$. Hence
$$
X^{\infty}_l \, = \, \bigcup_{\iota(\tau_p) \leq l}
W^{\infty}_u(\tau_p^{\infty})\,\cup \,\bigcup_{w <
\tau_p^{\infty}} W_u(w).
$$
It follows  that $X^{\infty}_l$ is a stratified set of top
dimension $ \leq l$. Without loss of generality, we may assume it
equal to $l$,
therefore $C(X^{\infty}_l)$ is  also a stratified set of top  dimension $l + 1$.\\
Now we use the gradient flow of $- \n J$ to deform
$C(X^{\infty}_l)$. By transversality arguments we can assume that
the deformation avoids all critical as well as critical points at
infinity having their Morse indices greater than $l + 2$. It
follows then, by  a Theorem of Bahri and Rabinowitz \cite{BR},
that  $C(X^{\infty}_l)$ retracts by deformation on the  set
\begin{equation}\label{f:fo}
U :=  X^{\infty}_l \,  \cup \, \bigcup_{morse(x{_{\infty}}) = l +
1} W^{\infty}_u(x_{\infty})  \, \cup \, \bigcup_{w <
\tau_p^{\infty}} W_u(w).
\end{equation}

\noindent Now taking $l = k - 1$ and using the fact that, by
assumption of theorem \ref{t:main}, there are no critical point at
 infinity with index $k$, we derive that $C(X^{\infty}_{k-1})$
retracts by deformation onto
\begin{equation}\label{f:rpdef}
Z_k^{\infty} := X^{\infty}_{k-1}  \, \cup  \bigcup_{w ; \,J'(w) =
0 ; \,w \,\mbox{\scriptsize dominated by } C(X^{\infty}_{k-1})  }
W_u(w).
\end{equation}
Now,  observe that it follows from the above deformation retract,
that the problem $(P_K)$ has necessary a solution $w$ with
$morse(w) \leq k$. Otherwise
 it follows from \eqref{f:rpdef} that

$$
1 \, = \chi(Z_k^{\infty}) \, = \,  \sum_{\tau_p \in \mathcal{F_1}
\, ;\, \iota(\tau_p) \leq  k - 1} (-1)^{ \iota(\tau_p)} \, ,
$$
where  $\chi$ denotes the Euler Characteristic. Such an equality contradicts the
second assumption of the theorem. \\
Now for generic $K$, it follows from the Sard-Smale Theorem that all solutions of
$(P_K)$ are nondegenerate solutions, in the sense
 that their associated linearized operator does not admit zero as an eigenvalue.
 See please \cite{SZ} for a related discussion in the Riemannian setting.\\
We derive now from \eqref{f:rpdef}, taking the Euler
Characteristic of both sides, that:
$$
1 \, = \chi(Z_k^{\infty}) \, = \,  \sum_{\tau_p \in \mathcal{F_1}
\, ;\, \iota(\tau_p) \leq  k - 1} (-1)^{ \iota(\tau_p)} \,  + \,
\sum_{w < C(X_{k-1}^{\infty}) ; \, J'(w) = 0} (-1)^{morse(w)}.
$$
It follows then that
$$
| 1 -\sum_{\tau_p \in \mathcal{F_1} \, ;\, \iota(\tau_p) \leq  k -
1} (-1)^{ \iota(\tau_p)} | \, \leq  \, | \sum_{w  ; \, J'(w) =
0,\, morse(w) \leq k} (-1)^{morse(w)}| \leq \#\mathcal{N}_k,
$$
where $\mathcal{N}_k $ denotes the set of solutions of $(P_K)$ having
their morse indices $\leq k$.
\end{pfn}

\section{A general existence result}

In this last section of this paper, we give a generalization of theorem \ref{t:main}.
 Namely instead of assuming that there are no   critical points at infinity
 of index $k$, we assume that the intersection number modulo 2, between the
 suspension of the complex at infinity of order $k$, $C(X_{k-1}^{\infty})$ and the
 stable manifold of all critical points at infinity of index $k $ is equal to zero.
 More precisely, for $\tau_p \in \mathcal{F_1}$ such that $\iota(\tau_p) \, = k $,
 we define the following intersection number:
$$
\mu_k(\tau_p) := C(X_{k - 1}^{\infty}) \, . \,
W_s^{\infty}(\tau_p^{\infty}) \qquad (\mbox{mod} 2).
$$
Observe that this intersection number is well defined since we may
assume by transversality that:
$$
\partial C(X_{k-1}^{\infty}) \, \cap \, W_s^{\infty}(\tau_p^{\infty}) \, = \, \emptyset.
$$
indeed $\,  dim (\partial C(X_{k-1}^{\infty})) \, = \,  k - 1$, {
while } \,$   codim (W_s^{\infty}(\tau_p^{\infty})) = k$.\\
 We are
now ready to state the following existence result:

\begin{thm}\label{t:a}

Let $0 < K \in C^2(M)$ be a  positive function satisfying the
conditions $(C_0)$ and $(C_1)$. If there exists $k \in \N$ such
that
\begin{enumerate}
\item
$$
\sum_{\tau_p \in \mathcal{F_1} ; \iota(\tau_p) \leq k - 1 } (-1)^{\iota(\tau_p)} \, \neq \, 1,
$$
\item
$$
\forall \tau_p \in \mathcal{F_1}, \mbox{ such that } \iota(\tau_p) = k,
\mbox{ there holds } \mu_k(\tau_p) \, = \, 0.
$$
\end{enumerate}
Then, there exists a solution $w$ of  the problem $(P_K)$ such
that:
$$
morse(w) \, \leq k.
$$
 Moreover, for generic
$K$ it holds
$$
\# \mathcal{N}_k \, \geq | 1 \, - \, \sum_{\tau_p \in
\mathcal{F_1} ;\, \iota(\tau_p) \leq k - 1 } (-1)^{\iota(\tau_p)}
|,
$$
where  $\mathcal{N}_k$ denotes the set of solutions of $(P_K)$
having their Morse indices less than or equal to $k$.
\end{thm}

\begin{pf}
The proof goes  along with the proof of theorem \ref{t:main},
therefore  we will only sketch the differences. Keeping the
notations of the proof of theorem \ref{t:main}, we observe that,
since

$$
\forall \tau_p \in \mathcal{F_1}, \mbox{ such that } \iota(\tau_p) = k, \mbox{ there holds } \mu_k(\tau_p) \, = \, 0,
$$
we may assume that the deformation of $C(X_{k-1}^{\infty})$ along
any pseudogradient flow of $-J$, avoids all critical points at
infinity
 having their Morse indices equal to $k$. It follows then from \eqref{f:fo}
that $C(X_{k-1}^{\infty})$
 retracts by deformation onto
\begin{equation}\label{f:rpdefa}
Z_k^{\infty} := X^{\infty}_{k-1}  \, \cup  \bigcup_{w ; \, J'(w) =
0 ; \,w \,\mbox{\scriptsize dominated by } C(X^{\infty}_{k-1})  }
W_u(w).
\end{equation}
Now the remainder of the proof is identical to the proof of theorem \ref{t:main}.
\end{pf}

\end{document}